# Reduced Differential Transform Method for (2+1) Dimensional type of the Zakharov–Kuznetsov ZK(*n,n*) Equations

Omer ACAN[1,a)] and Yıldıray KESKİN[1,b)]

[1)]*Department of Mathematics, Science Faculty, Selcuk University, Konya 42003, Turkey.*

[a)]acan_omer@ selcuk.edu.tr and [b)] yildiraykeskin@yahoo.com

**Abstract.** In this paper, reduced differential transform method (RDTM) is employed to approximate the solutions of (2+1) dimensional type of the Zakharov–Kuznetsov partial differential equations. We apply these method to two examples. Thus, we have obtained numerical solution partial differential equations of Zakharov–Kuznetsov. These examples are prepared to show the efficiency and simplicity of the method.

## 1. INTRODUCTION

Partial differential equations (PDEs) have numerous essential applications in various fields of science and engineering such as fluid mechanic, thermodynamic, heat transfer, physics [1]. Most new nonlinear PDEs do-not have a precise analytic solution. so numerical methods have largely been used to handle these equation. It is difficult to handle nonlinear part of these equations. Although most of scientists applied numerical methods to find the solution of these equations, solving such equations analytically is of fundamental importance since the existent numerical methods which approximate the solution of PDEs don't result in such an exact and analytical solution which is obtained by analytical methods.

In recent years, many researchers have paid attention to studying the solutions of nonlinear PDEs by various methods, for example, the Adomian's decomposition method (ADM) [2-7], the variational iteration method (VIM) [8-11], the homotopy analysis method [12-14], the homotopy perturbation method (HPM) [15-17], the differential transform method (DTM) [18-21], Hirtoa's bilinear method [22], the balance method [23], inverse scattering method [24], The RDTM was first proposed by Y. Keskin [25-28] in 2009. It has received much attention since it has applied to solve a wide variety of problems by many authors [29–33].

Zakharov–Kuznetsov (ZK) equation,

$$u_t + auu_x + (\nabla^2 u)_x = 0 \qquad (1.1)$$

where $\nabla^2 = \partial_x^2 + \partial_y^2$ or $\nabla^2 = \partial_x^2 + \partial_y^2 + \partial_z^2$ is the isotropic Laplacian [34-36]. This means that the ZK equation is given by

$$u_t + auu_x + (u_{xx} + u_{yy})_x = 0, \qquad (1.2)$$

and

$$u_t + auu_x + (u_{xx} + u_{yy} + u_{zz})_x = 0 \qquad (1.3)$$

in (2+1) and (3+1)-dimensional spaces. The ZK equation was first derived for describing weakly nonlinear ion-acoustic wave in a strongly magnetized lossless plasma in two dimensions [34] A further discussion of the analytical properties of the ZK equation and some constructive results were given.

Recently, in 2005, A.M. Wazwaz [37] studied the type of Zakharov–Kuznetsov equation, that is, the (2+1) dimensional and (3+1) dimensional ZK(*n,n*) equations of the form

$$u_t + a(u^n)_x + b(u^n)_{xxx} + k(u^n)_{yyx} = 0, \quad b,k > 0 \tag{1.4}$$

and

$$u_t + a(u^n)_x + b(u^n)_{xxx} + k(u^n)_{yyx} + r(u^n)_{zzx} = 0, \quad b,k,r > 0 \tag{1.5}$$

and in 2007, C. Lin, X. Zhang, [39] studied the (3+1) dimensional modified ZK equation of the form

$$u_t + au^p u_x + (u_{xx} + u_{yy} + u_{zz})_x = 0. \tag{1.6}$$

The main purpose of this study paper has been organized as follows: Section 2 deals with the analysis of the method. In Section 3, we apply the RDTM to solve types of ZK(*n,n*) equations of the form (1.5) equations. Conclusions are given in Section 4.

## 2. ANALYSIS OF THE RDTM

The basic definitions in the RDTM [25] are as follows:
**Definition 2.1.**

If function $u(x,t)$ is analytic and differentiated continuously with respect to time t and space x in the domain of interest, then let

$$U_k(x) = \frac{1}{k!}\left[\frac{\partial^k}{\partial t^k} u(x,t)\right]_{t=0} \tag{2.1}$$

where the t-dimensional spectrum function $U_k(x)$ is the transformed function. In this paper, the lowercase $u(x,t)$ represent the original function while the uppercase $U_k(x)$ stand for the transformed function.

The differential inverse transform of $U_k(x)$ is defined as follows:

$$u(x,t) = \sum_{k=0}^{\infty} U_k(x) t^k. \tag{2.2}$$

Then combining equation (2.1) and (2.2) we write

$$u(x,t) = \sum_{k=0}^{\infty} \frac{1}{k!}\left[\frac{\partial^k}{\partial t^k} u(x,t)\right]_{t=0} t^k. \tag{2.3}$$

From the above definitions, it can be found that the concept of the reduced differential transform is derived from the power series expansion.

For the purpose of illustration of the methodology to the proposed method, we write the ZK(*n,n*) equation in the standard operator form

$$L(u(x,t)) + N(u(x,t)) = 0 \tag{2.4}$$

with initial condition

$$u(x,y,0) = f(x,y) \tag{2.5}$$

where $L = \frac{\partial}{\partial t}$ is a linear operator which has partial derivatives, $N(u(x,t)) = a(u^n)_x + b(u^n)_{xxx} + k(u^n)_{yyx}$ is a nonlinear term.

**TABLE 1.** The fundamental operators of RDTM [25-28]

| Functional Form | Transformed Form |
|---|---|
| $u(x,y,t)$ | $U_k(x,y) = \frac{1}{k!}\left[\frac{\partial^k}{\partial t^k} u(x,y,t)\right]_{t=0}$ |
| $w(x,y,t) = u(x,y,t) \pm \alpha v(x,y,t)$ | $W_k(x,y) = U_k(x,y) \pm \alpha V_k(x,y)$ ($\alpha$ is a constant) |
| $w(x,y,t) = x^m y^n t^p u(x,y,t)$ | $W_k(x,y) = x^m y^n U_{(k-p)}(x,y)$ |
| $w(x,y,t) = u(x,y,t) v(x,y,t)$ | $W_k(x,y) = \sum_{r=0}^{k} V_r(x,y) U_{k-r}(x,y) = \sum_{r=0}^{k} U_r(x,y) V_{k-r}(x,y)$ |

$$w(x,y,t) = \frac{\partial^r}{\partial t^r} u(x,y,t) \qquad W_k(x,y) = (k+1)...(k+r)U_{k+r}(x,y) = \frac{(k+r)!}{k!}U_{k+r}(x,y)$$

$$w(x,y,t) = \frac{\partial}{\partial x^m \partial y^n} u(x,y,t) \qquad W_k(x,y) = \frac{\partial}{\partial x^m \partial y^n} U_k(x,y)$$

According to the RDTM and Table 1, we can construct the following iteration formula:

$$(k+1)U_{k+1}(x) = -N(U_k(x)) \tag{2.6}$$

where $N_k = N(U_k(x))$ is the transformations of the function $N(u(x,t))$ respectively.

For the easy to follow of the reader, we can give the first few nonlinear term are

$$N_0 = a\frac{\partial}{\partial x}U_0^n(x,y) + b\frac{\partial^3}{\partial x^3}U_0^n(x,y) + k\frac{\partial^3}{\partial y^2 \partial x}U_0^n(x,y)$$

$$N_1 = a\frac{\partial}{\partial x}\left(n\ U_0^{(n-1)}(x,y)U_1(x,y)\right) + b\frac{\partial^3}{\partial x^3}\left(n\ U_0^{(n-1)}(x,y)U_1(x,y)\right) + k\frac{\partial^3}{\partial y^2 \partial x}\left(n\ U_0^{(n-1)}(x,y)U_1(x,y)\right)$$

$$N_2 = a\frac{\partial}{\partial x}\left(\frac{n}{2}\left(n\ U_0^{(n-2)}(x,y)U_1^2(x,y) + n\ U_0^{(n-1)}(x,y)U_2(x,y) - U_0^{(n-2)}(x,y)U_1^2(x,y)\right)\right)$$

$$+ b\frac{\partial^3}{\partial x^3}\left(\frac{n}{2}\left(n\ U_0^{(n-2)}(x,y)U_1^2(x,y) + n\ U_0^{(n-1)}(x,y)U_2(x,y) - U_0^{(n-2)}(x,y)U_1^2(x,y)\right)\right)$$

$$+ k\frac{\partial^3}{\partial y^2 \partial x}\left(\frac{n}{2}\left(n\ U_0^{(n-2)}(x,y)U_1^2(x,y) + n\ U_0^{(n-1)}(x,y)U_2(x,y) - U_0^{(n-2)}(x,y)U_1^2(x,y)\right)\right)$$

**Maple Code for Nonlinear Function** as given [25]

```
restart;
NF:=Nu(t,x,y)^p:#Nonlinear Function
m:=4:        # Order
u[t]:=sum(u[b]*t^b,b=0..m):
NF[t]:=subs(Nu(x,t)=u[t],NF):
s:=expand(NF[t],t):
dt:=unapply(s,t):
for i from 0  to m do
n[i]:=((D@@i)(dt)(0)/i!):
print(N[i],n[i]); #Transform Function
od:
```

From initial condition (1.2), we write

$$U_0(x,y) = f(x,y) \tag{2.7}$$

Substituting (2.7) into (2.6) and by a straight forward iterative calculations, we get the following $U_k(x,y)$ values.

Then the inverse transformation of the set of values $\{U_k(x,y)\}_{k=0}^n$ gives approximation solution as,

$$\tilde{u}_n(x,y,t) = \sum_{k=0}^n U_k(x,y)t^k \tag{2.8}$$

where $n$ is order of approximation solution.

Therefore, the exact solution of problem is given by

$$u(x,y,t) = \lim_{n \to \infty} \tilde{u}_n(x,y,t). \tag{2.9}$$

## 3. NUMERICAL APPLICATIONS

In this section, we test the RDTM for the ZK(*3,3*) and ZK(*2,2*) equations whit fully nonlinear dispersion. Numerical results are very encouraging.

**Example 3.1.** First we consider the following ZK(*3,3*) equations [38]:

$$u_t + (u^3)_x + 2\ u^3)_{xxx} + 2\ u^3)_{yyx} = 0 \tag{3.1}$$

subject to initial condition

$$u(x,y,0) = \frac{3}{2}\lambda \sinh\frac{1}{6}(x+y) \tag{3.2}$$

Taking differential transform of (3.1) and the initial condition (3.2) respectively, we obtain

$$U_{k+1}(x,y) = -\frac{1}{(k+1)}\left(\frac{\partial}{\partial x}U_0^3(x) + 2\frac{\partial^3}{\partial x^3}U_0^3(x) + 2\frac{\partial^3}{\partial y^2 \partial x}U_0^3(x)\right) \tag{3.3}$$

where the t-dimensional spectrum function $U_k(x,y)$ are the transformed function.

From the initial condition (3.2) we write

$$U_0(x,y) = \frac{3}{2}\lambda \sinh\frac{1}{6}(x+y). \tag{3.4}$$

Now, substituting (3.4) into (3.3), we obtain the following $U_k(x,y)$ values successively

$$U_1(x,y) = -\frac{3}{8}\lambda^3 \cosh\left(\frac{x+y}{6}\right)\left(9\cosh^3\left(\frac{x+y}{6}\right) - 8\right)$$

$$U_2(x,y) = \frac{3}{64}\lambda^2 \sinh\left(\frac{x+y}{6}\right)\left(765\cosh^4\left(\frac{x+y}{6}\right) - 729\cosh^2\left(\frac{x+y}{6}\right) + 91\right)$$

$$U_3(x,y) = -\frac{1}{256}\lambda^7 \cosh\left(\frac{x+y}{6}\right)\begin{pmatrix} -382293\cosh^4\left(\frac{x+y}{6}\right) + 188181\cosh^6\left(\frac{x+y}{6}\right) \\ +234468\cosh^2\left(\frac{x+y}{6}\right) - 39851 \end{pmatrix}$$

$$U_4(x,y) = \frac{1}{4096}\lambda^9 \sinh\left(\frac{x+y}{6}\right)\begin{pmatrix} 93534345\cosh^8\left(\frac{x+y}{6}\right) - 198626022\cosh^6\left(\frac{x+y}{6}\right) \\ +135212355\cosh^4\left(\frac{x+y}{6}\right) - 30715929\cosh^2\left(\frac{x+y}{6}\right) + 1179946 \end{pmatrix}$$

(3.5)

From (2.8)

$$\tilde{u}_4(x,y,t) = \sum_{k=0}^{4} U_k(x,y)t^k \tag{3.6}$$

Substituting (3.4) and (3.5) in (3.6), we have

$$\tilde{u}_4(x,y,t) = \frac{1}{4096}\lambda \begin{pmatrix} 6144\lambda^2 t \sinh\left(\frac{x+y}{6}\right) - 13824\cosh^3\left(\frac{x+y}{6}\right) + 12288\lambda^2 t \cosh\left(\frac{x+y}{6}\right) \\ +146880\lambda^4 t^2 \sinh\left(\frac{x+y}{6}\right)\cosh^4\left(\frac{x+y}{6}\right) - 139968\lambda^4 t^2 \sinh\left(\frac{x+y}{6}\right)\cosh^2\left(\frac{x+y}{6}\right) \\ +17472\lambda^4 t^2 \sinh\left(\frac{x+y}{6}\right) + 6116688\lambda^6 t^3 \cosh^5\left(\frac{x+y}{6}\right) - 3010896\lambda^6 t^3 \cosh^7\left(\frac{x+y}{6}\right) \\ -3751488\lambda^6 t^3 \cosh^3\left(\frac{x+y}{6}\right) + 637616\lambda^9 t^3 \cosh\left(\frac{x+y}{6}\right) + 93534345\lambda^8 t^4 \sinh\left(\frac{x+y}{6}\right)\cosh^8\left(\frac{x+y}{6}\right) \\ -198626022\lambda^8 t^4 \sinh\left(\frac{x+y}{6}\right)\cosh^6\left(\frac{x+y}{6}\right) + 135212355\lambda^8 t^4 \sinh\left(\frac{x+y}{6}\right)\cosh^4\left(\frac{x+y}{6}\right) \\ -30715929\lambda^8 t^4 \sinh\left(\frac{x+y}{6}\right)\cosh^2\left(\frac{x+y}{6}\right) + 1179946\lambda^8 t^4 \sinh\left(\frac{x+y}{6}\right) \end{pmatrix}$$

when analyzed the solutions of ZK(*3,3*) equations by RDTM, the following results are obtained:

**TABLE 1.** For ZK(*3,3*), comparison of absolute error of the numerical results for $\tilde{u}_4(x,y,t)$, by RDTM

| $\lambda$ | $x$ | $y$ | $t$ | RDTM Solution | Abs-Error-RDTM |
|---|---|---|---|---|---|
| | 0.0 | 0.0 | | $-0.375000000 \times 10^{-18}$ | $0.2265468750 \times 10^{-37}$ |
| | 0.0 | 0.5 | | $0.1251447262 \times 10^{-5}$ | $0.2292047232 \times 10^{-37}$ |
| | 0.0 | 1.0 | | $0.2511590160 \times 10^{-5}$ | $0.2373142696 \times 10^{-37}$ |
| | 0.5 | 0.0 | | $0.1251447262 \times 10^{-5}$ | $0.2292047232 \times 10^{-37}$ |
| 0.00001 | 0.5 | 0.5 | 0.001 | $0.2511590160 \times 10^{-5}$ | $0.2373142096 \times 10^{-37}$ |
| | 0.5 | 1.0 | | $0.3789184752 \times 10^{-5}$ | $0.2512918476 \times 10^{-37}$ |
| | 1.0 | 0.0 | | $0.2511590160 \times 10^{-5}$ | $0.2373142696 \times 10^{-37}$ |
| | 1.0 | 0.5 | | $0.3789184752 \times 10^{-5}$ | $0.2512918476 \times 10^{-37}$ |
| | 1.0 | 1.0 | | $0.5093108360 \times 10^{-5}$ | $0.2718596318 \times 10^{-37}$ |

**Example 3.2.** Now we consider the following ZK(*2,2*) equation:

$$u_t + (u^2)_x + \frac{1}{8}(u^2)_{xxx} + \frac{1}{8}(u^2)_{yyx} = 0 \tag{3.7}$$

subject to initial condition

$$u(x,y,0) = -\frac{4}{3}\lambda \cosh^2(x+y) \tag{3.8}$$

Similarly by using the RDTM and Table 1 to (3.7) and (3.8), we obtain the recursive relation

$$U_{k+1}(x,y) = -\frac{1}{(k+1)}\left(\frac{\partial}{\partial x}U_0^2(x) + \frac{1}{8}\frac{\partial^3}{\partial x^3}U_0^2(x) + \frac{1}{8}\frac{\partial^3}{\partial y^2 \partial x}U_0^2(x)\right) \tag{3.9}$$

where the t-dimensional spectrum function $U_k(x,y)$ are the transformed function.

From the initial condition (3.8) we write

$$U_0(x,y) = -\frac{4}{3}\lambda \cosh^2(x+y) \tag{3.10}$$

Now, substituting (3.10) into (3.9), we obtain the following $U_k(x,y)$ values successively

$$U_1(x,y) = -\frac{32}{9}\lambda^2 \cosh(x+y)\sinh(x+y)\left(10\cosh^2(x+y)-3\right)$$

$$U_2(x,y) = -\frac{64}{27}\lambda^3 \left(1200\cosh^6(x+y)-1520\cosh^4(x+y)+408\cosh^2(x+y)-9\right)$$

$$U_3(x,y) = -\frac{4096}{243}\lambda^4 \cosh(x+y)\sinh(x+y)\begin{pmatrix}23800\cosh^6(x+y)-28500\cosh^4(x+y)\\+8265\cosh^2(x+y)-423\end{pmatrix}$$

$$U_4(x,y) = -\frac{1024}{729}\lambda^5 \begin{pmatrix}-11151+1512792\cosh^2(x+y)-124257600\cosh^8(x+y)+58864000\cosh^{10}(x+y)\\-21520320\cosh^4(x+y)+85809600\cosh^6(x+y)\end{pmatrix}$$

(3.11)

From (2.8)

$$\tilde{u}_4(x,y,t) = \sum_{k=0}^{4} U_k(x,y)t^k \tag{3.12}$$

Substituting (3.10) and (3.11) in (3.12), we have

$$\tilde{u}_4(x,y,t) = \frac{4}{729}\lambda \begin{pmatrix} 243\cosh^2(x+y) + 6480\lambda t\cosh^3(x+y)\sinh(x+y) - 656640\lambda^2 t^2\cosh^4(x+y) \\ +176256\lambda^2 t^2\cosh^2(x+y) + 3888\lambda^2 t^2 + 73113600\lambda^3 t^3\cosh^7(x+y)\sinh(x+y) \\ -87552000\lambda^3 t^3\cosh^5(x+y)\sinh(x+y) + 25390080\lambda^3 t^3\cosh^3(x+y)\sinh(x+y) \\ -1299456\lambda^3 t^3\cosh(x+y)\sinh(x+y) - 2854656\lambda^4 t^4 + 387274752\lambda^4 t^4 c\,osh^2(x+y) \\ -31809945600\lambda^4 t^4 c\,osh^8(x+y) + 15069184000\lambda^4 t^4 c\,osh^{10}(x+y) \\ -5509201920\lambda^4 t^4 c\,osh^4(x+y) + 21967257600\lambda^4 t^4 c\,osh^6(x+y) \end{pmatrix}$$

when analyzed the solutions of ZK(2,2) equations by RDTM, the following results are obtained:

**TABLE 2.** For ZK(2,2), comparison of absolute error of the numerical results for $\tilde{u}_4(x,y,t)$, by RDTM

| $\lambda$ | $x$ | $y$ | $t$ | RDTM Solution | Abs-Error-RDTM |
|---|---|---|---|---|---|
|  | 0.0 | 0.0 |  | -0.00001333333333 | $0.2932786015 \times 10^{-24}$ |
|  | 0.0 | 0.5 |  | -0.00001695387292 | $0.2692669148 \times 10^{-20}$ |
|  | 0.0 | 1.0 |  | -0.00003174798469 | $0.1562939098 \times 10^{-19}$ |
|  | 0.5 | 0.0 |  | -0.00001695387292 | $0.2692669148 \times 10^{-20}$ |
| 0.00001 | 0.5 | 0.5 | 0.001 | -0.00003174798469 | $0.1562939098 \times 10^{-19}$ |
|  | 0.5 | 1.0 |  | -0.00007378450649 | $0.1005498720 \times 10^{-18}$ |
|  | 1.0 | 0.0 |  | -0.00003174798469 | $0.1562939098 \times 10^{-19}$ |
|  | 1.0 | 0.5 |  | -0.00007378450649 | $0.1005498720 \times 10^{-18}$ |
|  | 1.0 | 1.0 |  | -0.00001887222243 | $0.6999706576 \times 10^{-18}$ |

Now for a better understanding, two examples given above are as follows absolute error graphics:

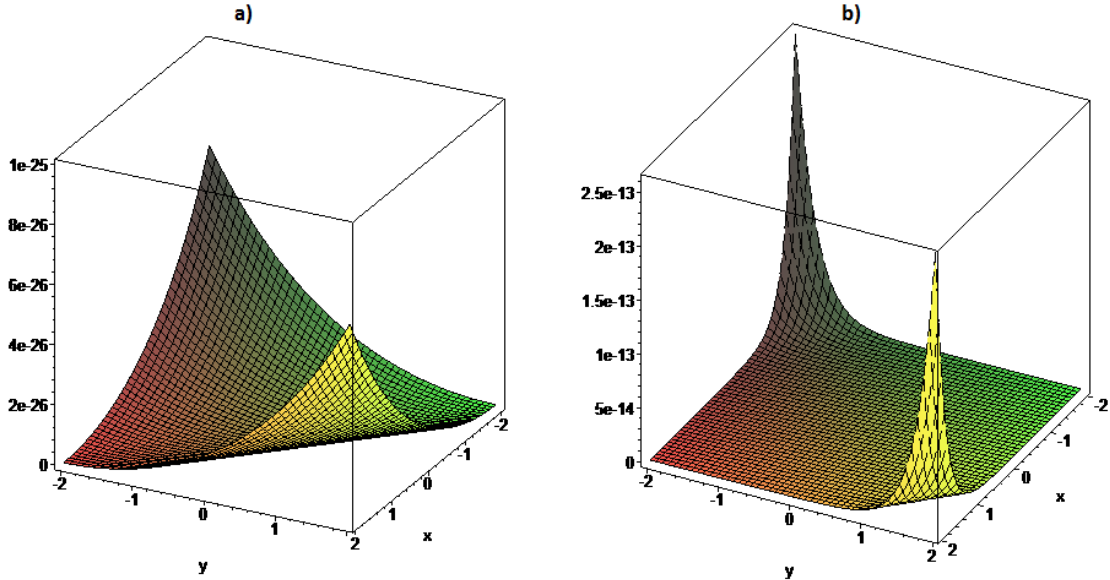

**FIG. 1:** Absolute errors for **(a)** ZK(3,3) and **(b)** ZK(2,2) equations when $-2 \leq x \leq 2$, $-2 \leq y \leq 2$, $t = 0,001$ and $\lambda = 0.00001$.

## 4. CONCLUSION

In this study, we apply RDTM on (2+1) dimensional ZK(3,3) and ZK(2,2) PDEs with fully nonlinear dispersion. The obtained result which obtained are highly reliable an encouraging.